# ON CONSTRUCTION OF GLOBAL ACTIONS OF FINITE PARTIAL GROUP ACTIONS ON SETS

Ram Parkash Sharma and Meenakshi
Department of Mathematics
Himachal Pradesh University Summerhill, Shimla-171005, India
E-mail: rp_math_hpu@yahoo.com

**Abstract**
A generalization of $G$–sets, called partial $G$–sets, are the sets that admit an action of partial maps on their subsets. Partial actions are a powerful tool to generalize many results of group actions. These generalizations are obtained by using global actions when they exist. The main objective of this paper is to construct the global action of a given finite partial group action on a set. For this, first we generalize the orbit-stabilizer theorem for partial group actions and use it to know the exact size of the orbits in the global set.

## 1. Introduction

First, we recall the definition of a partial action on a set from [3].

**Definition 1.1.** Let $G$ be a group and $X$, a set. A partial action of $G$ on $X$ is a pair $\alpha = \{\{D_g\}_{g \in G}, \{\alpha_g\}_{g \in G}\}$, where for each $g \in G$, $D_g$ is a subset of $X$ and $\alpha_g : D_{g^{-1}} \to D_g$ is a bijective map, satisfying the following three properties for each $g, h \in G$ :
(i) $D_1 = X$ and $\alpha_1 = Id_X$, the identity map on $X$,
(ii) $\alpha_g(D_{g^{-1}} \cap D_h) = D_g \cap D_{gh}$,
(iii) $\alpha_g(\alpha_h(x)) = \alpha_{gh}(x)$ for $x \in D_{h^{-1}} \cap D_{h^{-1}g^{-1}}$.

If $\alpha$ is a partial action of $G$ on $X$, then we say that $X$ is a partial $G$–set or $(X, \alpha)$ as a partial action. Partial actions of groups appeared in many areas of mathematics and are studied by many authors on different structures (see [3,4,5,6,7 and 8]). When a partial action on some structure is given, then one of the most relevant problem is the question of existence and uniqueness of a globalization; that is, a global group action (enveloping action) whose restriction to the original object gives initial partial action. In [1], Abadie discussed the problem of deciding whether a given partial action is the restriction of a global action and uniqueness of this global action. The definition of the global action of a partial action and the exact solution of the problem depend upon the category under consideration. For example, in the category of $C^*$–algebras, enveloping action is unique when it exists; whereas in the category of the topological spaces, for each partial action there always exists a unique enveloping action. In particular, when $X$ a is non-empty set and $\alpha$ is a partial action of $G$ on $X$, there exists a global action $(T,\beta)$ for $(X,\alpha)$, where $T$ is a set with global action $\beta = \{\beta_g \mid g \in G\}$ on it such that $\beta_g's$ are bijections on $T$ and $\alpha$ is the restriction of $\beta$ on $X$. That is, $X \subseteq T$, $T$ is the orbit of $X$, $D_g = X \cap \beta_g(X)$ for each $g \in G$ and $\alpha_g(x) = \beta_g(x) \forall x \in D_{g^{-1}}$. In this situation, for each $x \in X$, the partial orbit $O_x^\alpha = \{\alpha_{g^{-1}}(x) \mid x \in D_g\} = X \cap O_x$ and the partial stabilizer $G_x^\alpha = \{g \in G \mid x \in D_{g^{-1}}$ and $\alpha_g(x) = x\}$ $= G_x$ as observed in [2]. Throughout this paper, the group $G$ is finite and for $x \in X$ we use $G_x$ for both $G_x^\alpha$ and $G_x$ as they are equal.

Many properties of partial actions are obtained by using their relationships to global actions, so global actions are considered important for partial actions when they exist. It is pertinent to note that the orbit-stabilizer theorem provides the order of an orbit $O_x$ in terms of the stabilizer $G_x$ and $G_x = G_x^\alpha$ when $x \in X$. Thus the orbit-stabilizer theorem of a partial action $(X,\alpha)$ can give information regarding the order of an orbit in $T$. But the theorem is not true for partial actions on sets as observed in Example 2.3. So, first we generalize the orbit-stabilizer theorem for partial actions on sets (Theorem 2.10). Then using this, we get $|O_x| = |O_x^\alpha| + |\overline{G^x}|/|G_x|$, where $\overline{G^x} = \{h \in G \mid x \notin D_{h^{-1}}\}$ (see Corollary 2.12). Finally, with the known size of the orbit $O_x$, the global action of a given partial action on a set is constructed.

## 2. Global Actions of Partial Actions on Sets

In this section, we consider a partial action $(X,\alpha)$ of a finite group $G$ on set $X'$ with the global action $(T,\beta)$.

**Definition 2.1.** Let $X'$ be a subset of partial $G$-set $X$. Then $X'$ is said to be a partial $G$-subset of $X$ if $x \in X' \cap D_g$ implies $\alpha_{g^{-1}}(x) \in X'$.

A natural example of a partial $G$-set is provided by the following result.

**Lemma 2.2.** For each $x \in X$, $O_x^\alpha$ is a partial $G$-subset of $X$ and hence $G$ acts partially on $O_x^\alpha$.

**Proof.** Let $x \in O_x^\alpha \cap D_g$. Then $\beta_{g^{-1}}(x) \in O_x$, as $O_x$ is a $G$-subset of $T$. Hence, $\alpha_{g^{-1}}(x) = \beta_{g^{-1}}(x) \in X \cap O_x = O_x^\alpha$.

The orbit-stabilizer theorem is not true, in general, for partial actions as observed below:
**Example 2.3.(i)** Let $X = \{x_1, x_2, x_3, x_4\}$ be a set. Let $G = \{g \mid g^8 = 1\}$ be a cyclic group of order $8, I = \{x_1, x_2, x_4\}, J = \{x_1, x_2, x_3\}, M = \{x_1, x_2\}$. Consider a partial action $\alpha$ of $G$ on $X$ given by $D_1 = X, D_g = D_{g^3} = D_{g^5} = D_{g^7} = \phi; D_{g^2} = I, D_{g^4} = M, D_{g^6} = J$ and
$\alpha_1 = Id_X; \alpha_g = \alpha_{g^3} = \alpha_{g^5} = \alpha_{g^7} = \phi; \alpha_{g^2} : x_1 \mapsto x_2 \mapsto x_1, x_3 \mapsto x_4; \alpha_{g^4} : x_1 \mapsto x_1, x_2 \mapsto x_2;$
$\alpha_{g^6} : x_1 \mapsto x_2 \mapsto x_1, x_4 \mapsto x_3$. We have $G_{x_1} = \{1, g^4\}$ so $[G : G_{x_1}] = 4$, but $O_{x_1}^\alpha = \{x_1, x_2\}$. Similarly, $G_{x_3} = \{1\}$ and $O_{x_3}^\alpha = \{x_3, x_4\}$. Hence, no element of $X$ satisfies the orbit-stabilizer theorem.

**(ii)** Let $X, I, J, M$ be as defined above and $G = \{g \mid g^4 = 1\}$ be a cyclic group of order four. Let $D_g = I, D_{g^2} = M$ and $D_{g^3} = J$ with $\alpha_g : x_1 \mapsto x_2 \mapsto x_1, x_3 \mapsto x_4;$
$\alpha_{g^2} : x_1 \mapsto x_1, x_2 \mapsto x_2; \alpha_{g^3} : x_1 \mapsto x_2 \mapsto x_1, x_4 \mapsto x_3$.
Here, $O_{x_1}^\alpha = \{x_1, x_2\} = O_{x_2}^\alpha, G_{x_1} = G_{x_2} = \{1, g^2\}$, then we have
$$G/G_{x_1} = \{G_{x_1}, g\, G_{x_1}\} \cong O_{x_1}^\alpha$$

under the map $x_1 \mapsto G_{x_1}$ and $x_2 \mapsto g\, G_{x_1}$. Further, $O_{x_3}^\alpha = \{x_3, x_4\}$ and $G_{x_3} = \{1\}$ and hence $[G : G_{x_3}] = 4$.

These examples suggest that alike form of the orbit-stabilizer theorem may exist for the partial actions on sets if for each $x \in X$, we replace $G$ by a subset $G^x$ of $G$ such that $x \in \bigcap_{g^{-1} \in G^x} D_g$. Therefore, we set $G^x = \{g^{-1} \mid x \in D_g\}$ and $G^x/G_x = \{g^{-1}G_x \mid x \in D_g\}$.

We observe that there is a natural partial action of $G$ on the set $G^x/G_x$.

**Lemma 2.4.** Let $x \in X$ and $G^x/G_x$ be as defined above. Then $G^x/G_x$ is a partial $G$–set.

**Proof.** For $x \in X$, there exists a natural group action $\overline{\beta}$ of group $G$ on $G/G_x$, that is, $\overline{\beta} : G/G_x \to G/G_x$ given by $\overline{\beta}_h(g^{-1}G_x) = hg^{-1}G_x$. The group action $\overline{\beta}$ induces a partial group action $\overline{\alpha}$ on $G^x/G_x$ by restricting $\overline{\beta}$ on $G^x/G_x$, that is,
$$\overline{D}_{h^{-1}} = G^x/G_x \cap \overline{\beta}_{h^{-1}}(G^x/G_x)$$
and the map $\overline{\alpha} : \overline{D}_{h^{-1}} \to \overline{D}_h$ is defined by $\overline{\alpha}_h(g^{-1}G_x) = hg^{-1}G_x, g^{-1}G_x \in \overline{D}_{h^{-1}}$.

The following result establishes a relation between the partial actions $\alpha$ and $\overline{\alpha}$.

**Lemma 2.5.** Let $X$ be a partial $G$–set and for $x \in X, (G^x/G_x, \overline{\alpha})$ be the partial action as given above. Then $\overline{D}_{h^{-1}} = \{g^{-1}G_x \mid x \in D_g \cap D_{gh^{-1}}\}$.

**Proof.** For the inclusion $\{g^{-1}G_x \mid x \in D_g \cap D_{gh^{-1}}\} \subseteq \overline{D}_{h^{-1}}$, we consider an element $g^{-1}G_x$ from the former. Then $x \in D_g \cap D_{gh^{-1}}$. Obviously, $g^{-1}G_x \in G^x/G_x$ as $x \in D_g$. By the definition of $\overline{\beta}$,

$$\overline{\beta}_{h^{-1}}(hg^{-1}G_x) = g^{-1}G_x.$$

Moreover, $x \in D_{gh^{-1}}$ implies that $hg^{-1}G_x \in G^x/G_x$, so that $g^{-1}G_x \in \overline{\beta}_{h^{-1}}(G^x/G_x) \cap G^x/G_x$.

On the other hand, let $g^{-1}G_x \in G^x/G_x \cap \overline{\beta}_{h^{-1}}(G^x/G_x)$. By the definition of $G^x/G_x$, we have $x \in D_g$. Further, $g^{-1}G_x \in \overline{\beta}_{h^{-1}}(G^x/G_x)$ implies that $\exists \ tG_x \in G^x/G_x$ such that

$$\overline{\beta}_{h^{-1}}(tG_x) = g^{-1}G_x.$$

Also,

$$\overline{\beta}_{h^{-1}}(hg^{-1}G_x) = g^{-1}G_x.$$

Since $\overline{\beta}$ is a global isomorphism, therefore we have

$$hg^{-1}G_x = tG_x \in G^x/G_x,$$

which implies that $x \in D_{gh^{-1}}$, proving the other inclusion.

**Remark 2.6.** Note that $g^{-1}G_x \in \overline{D}_{h^{-1}}$ implies that $x \in D_g \cap D_{gh^{-1}}$. Now $\overline{\alpha}_h(g^{-1}G_x) = g^{-1}G_x \in \overline{D}_{h^{-1}}$ if and only if $x \in D_{gh^{-1}} \cap D_{(hg^{-1})^{-1}h} = D_g \cap D_{gh^{-1}}$ which is true. Therefore the above result confirms that if $g^{-1}G_x \in \overline{D}_{h^{-1}}$, then $\overline{\alpha}_h(g^{-1}G_x) \in \overline{D}_h$.

**Definition 2.7.** Let $(X, \alpha)$ and $(X', \alpha')$ be two partial actions of a group $G$. Then a map $\Phi : X \to X'$ is said to be a partial $G$–map if $x \in D_{g^{-1}}$ implies that $\Phi(x) \in D'_{g^{-1}}$ and $\Phi(\alpha_g(x)) = \alpha'_g(\Phi(x))$. Further, $X$ is said to be partially $G$–isomorphic to $X'$ if $\Phi$ is a bijection.

**Remark 2.8.** A partial $G$–map can be obtained by restricting $G$– maps. Suppose that $\psi$ is a $G$– map between two global actions $(T, \beta)$ and $(T', \beta')$. Let $(X, \alpha)$ and $(X', \alpha')$ where $X \subseteq T, X' \subseteq T'$ with $\alpha$ and $\alpha'$ the restrictions of $\beta$ and $\beta'$ respectively be two partial actions. Let $\phi$ be the restriction of $\psi$ on $X$ such that $\phi(X) \subseteq X'$. Then $\phi$ is a partial action.

**Definition 2.9.** Let $(X,\alpha)$ be a partial action. Then a subset $S$ of $X$ is said to be a partial $G$-transversal in $X$ if it meets each partial orbit in $X$ exactly once.

Now we prove the orbit–stabilizer theorem for partial actions.

**Theorem 2.10 (The Orbit-Stabilizer Theorem).** Let $(X,\alpha)$ be a partial action. Then for $x \in X$, $O_x^\alpha$ is partially $G$-isomorphic to $G^x/G_x$ and hence

$$|O_x^\alpha| = |G^x|/|G_x|.$$

Further,

$$X \cong \bigcup_{s \in S} G^s/G_s,$$

where $S$ is a partial $G$-transversal in $X$.

**Proof.** Let $\Phi : O_x \to G/G_x$ be the $G$-isomorphism defined by

$$\Phi(\beta_{g^{-1}}(x)) = g^{-1}G_x.$$

For $y \in O_x^\alpha$, there exists $g \in G$ such that $y = \alpha_{g^{-1}}(x)$ and $x \in D_g$ so that $\Phi(\beta_{g^{-1}}(x)) \in G^x/G_x$.

Let $\phi$ be the restriction of $\Phi$ to $O_x^\alpha$, that is, $\phi : O_x^\alpha \to G^x/G_x$ be given by

$$\phi(\alpha_{g^{-1}}(x)) = g^{-1}G_x, x \in D_g.$$

As $\phi(O_x^\alpha) \subseteq G^x/G_x$, by Remark 2.8, $\phi$ is a partial $G$-map which is one-one.

In order to prove that $\phi$ is onto, it is sufficient to show that $\Phi^{-1}(G^x/G_x) \subseteq O_x^\alpha$. For, let $y \in O_x$ such that $\Phi(y) \in G^x/G_x$. Since $y \in O_x, \exists g \in G$ such that $y = \beta_{g^{-1}}(x)$ so $g^{-1}G_x \in G^x/G_x$. But this implies that $x \in D_g$ and $y = \alpha_{g^{-1}}(x) \in O_x^\alpha$. As orbits of the global action $(T,\beta)$ induce a partition of $T$, the partial orbits of $(X,\alpha)$ induce a partition of $X$. Hence, we get $X = \bigcup_{s \in S} O_x^\alpha \cong \bigcup_{s \in S} G^s/G_s$.

**Example 2.11.** Consider the partial action $\alpha$ of $G$ on $X$ given in Example 2.3(i). Here $O_{x_1}^\alpha = \{x_1, x_2\}$ and $O_{x_3}^\alpha = \{x_3, x_4\}$. Let $S = \{x_1, x_3\}$ be a partial $G$-transversal in $X$. Then $G^{x_1} = \{1, g^2, g^4, g^6\}$ and $G_{x_1} = \{1, g^4\}$. Since $g^4 G_{x_1} = G_{x_1}$ and $g^6 G_{x_1} = g^2 G_{x_1}$, we have $G^{x_1}/G_{x_1} = \{G_{x_1}, g^2 G_{x_1}\} \cong O_{x_1}^\alpha$ under the map $x_1 \mapsto G_{x_1}$ and $x_2 \mapsto g^2 G_{x_1}$. Similarly, $G^{x_3} = \{1, g^2\}, G_{x_3} = \{1\}$ and $G^{x_3}/G_{x_3} = \{G_{x_3}, g^2 G_{x_3}\} \cong O_{x_3}^\alpha$ under the map

$x_3 \mapsto G_{x_3}$ and $x_4 \mapsto g^2 G_{x_3}$. Hence, $X \cong \bigcup_{s \in S} G^s / G_s$ under the maps given above with $S = \{x_1, x_3\}$.

As an application of the orbit-stabilizer theorem, we get

**Corollary 2.12**. For $x \in X$, we have $|O_x| = |O_x^\alpha| + |\overline{G^x}|/|G_x|$, where $\overline{G^x} = \{h \in G \mid x \notin D_{h^{-1}}\}$.

**Proof.** By the orbit-stabilizer theorems for global actions and partial actions, we have $|G_x| \mid |G|$ and $|G_x| \mid |G^x|$ which imply that $|G_x| \mid (|G| - |G^x|) = |\overline{G^x}|$. Hence

$$|O_x| = |G|/|G_x|$$
$$= |G^x| + |\overline{G^x}|/|G_x|$$
$$= |G^x|/|G_x| + |\overline{G^x}|/|G_x|$$
$$= |O_x^\alpha| + |\overline{G^x}|/|G_x|.$$

We also need

**Lemma 2.13.** Let $(X, \alpha)$ be a partial action of a finite group $G$ with the global action $(T, \beta)$. Then there is a one-one correspondence between the $G$-orbits in $T$ and partial $G$-orbits in $X$. In particular, if $X$ is finite, then the number $k$ of disjoint orbits in $X$ is also given by:

$$k = \frac{1}{|G|} \sum_{g \in G} T_g, \quad T_g = \{u \in T \mid \beta_g(u) = u\}.$$

**Proof.** Let $O_t, t \in T$ be an orbit in $T$. Then for $t \in T, \exists g \in G$ such that $x = \beta_g(t) \in T$. Hence $x \in O_t$ so that $O_t = O_x, x \in X$. This sets up a one-one correspondence, via the map $O_t \to O_x \cap X = O_x^\alpha, x \in O_t \cap X$ between $G$-orbits in $T$ and partial $G$-orbits in $X$. The rest of the proof follows from Burnside Theorem.

Let $(X, \alpha)$ be a partial action of a finite group $G$ on a set $X$ and $(T, \beta)$ be its enveloping action. Let $S = \{x_1, x_2, \ldots x_m, \ldots\}$ be a partial $G$-transversal in $X$. Then $S$ is also a $G$-transversal in $T$, as observed in above lemma. Hence $O_{x_1}, O_{x_2}, \ldots O_{x_m}, \ldots$ are disjoint orbits in $T$ with the size given by Corollary 2.12 such that $T = \bigcup_i O_{x_i}$. Now, the complete action of $\beta$ on $T$ can be obtained using the orbits $O_{x_1}, O_{x_2}, \ldots O_{x_m}, \ldots$ and the partial action $(X, \alpha)$.

**Example 2.14**. Consider the partial action given in Example 2.3(1). We have $|O_{x_1}^\alpha| = 2, |G_{x_1}| = 2$ and $|\overline{G^{x_1}}| = 4$. By using Corollary 2.12, we have $|O_{x_1}| = 2 + 2 = 4$.

Similarly, $|O^\alpha_{x_3}|=2, |G_{x_3}|=1$ and $|\overline{G^{x_3}}|=6$. Therefore, $|O_{x_3}|=2+6=8$. The global action $(T,\beta)$ of the partial action $(X,\alpha)$ is given by:

The global set $T=\{x_1,x_2,x_3,x_4,x_5,x_6,x_7,x_8,x_9,x_{10},x_{11},x_{12}\}$ has two disjoint orbits $O_{x_1}=\{x_1,x_2,x_5,x_6\}$ and $O_{x_3}=\{x_3,x_4,x_7,x_8,x_9,x_{10},x_{11},x_{12}\}$. The bijections $\beta_g's$ are computed using $O_{x_1}, O_{x_3}$ and $\alpha_g's$ as follows:

| $\beta_g$ | $\beta_{g^2}$ | $\beta_{g^3}$ | $\beta_{g^4}$ | $\beta_{g^5}$ | $\beta_{g^6}$ | $\beta_{g^7}$ | $\beta_{g^8}=\beta_1$ |
|---|---|---|---|---|---|---|---|
| $x_1 \mapsto x_5$ | $x_1 \mapsto x_2$ | $x_1 \mapsto x_6$ | $x_1 \mapsto x_1$ | $x_1 \mapsto x_5$ | $x_1 \mapsto x_2$ | $x_1 \mapsto x_6$ | $x_1 \mapsto x_1$ |
| $x_2 \mapsto x_6$ | $x_2 \mapsto x_1$ | $x_2 \mapsto x_5$ | $x_2 \mapsto x_2$ | $x_2 \mapsto x_6$ | $x_2 \mapsto x_1$ | $x_2 \mapsto x_5$ | $x_2 \mapsto x_2$ |
| $x_3 \mapsto x_7$ | $x_3 \mapsto x_4$ | $x_3 \mapsto x_8$ | $x_3 \mapsto x_9$ | $x_3 \mapsto x_{10}$ | $x_3 \mapsto x_{11}$ | $x_3 \mapsto x_{12}$ | $x_3 \mapsto x_3$ |
| $x_4 \mapsto x_8$ | $x_4 \mapsto x_9$ | $x_4 \mapsto x_{10}$ | $x_4 \mapsto x_{11}$ | $x_4 \mapsto x_{12}$ | $x_4 \mapsto x_3$ | $x_4 \mapsto x_7$ | $x_4 \mapsto x_4$ |
| $x_5 \mapsto x_2$ | $x_5 \mapsto x_6$ | $x_5 \mapsto x_1$ | $x_5 \mapsto x_5$ | $x_5 \mapsto x_2$ | $x_5 \mapsto x_6$ | $x_5 \mapsto x_1$ | $x_5 \mapsto x_5$ |
| $x_6 \mapsto x_1$ | $x_6 \mapsto x_5$ | $x_6 \mapsto x_2$ | $x_6 \mapsto x_6$ | $x_6 \mapsto x_1$ | $x_6 \mapsto x_5$ | $x_6 \mapsto x_2$ | $x_6 \mapsto x_6$ |
| $x_7 \mapsto x_4$ | $x_7 \mapsto x_8$ | $x_7 \mapsto x_9$ | $x_7 \mapsto x_{10}$ | $x_7 \mapsto x_{11}$ | $x_7 \mapsto x_{12}$ | $x_7 \mapsto x_3$ | $x_7 \mapsto x_7$ |
| $x_8 \mapsto x_9$ | $x_8 \mapsto x_{10}$ | $x_8 \mapsto x_{11}$ | $x_8 \mapsto x_{12}$ | $x_8 \mapsto x_3$ | $x_8 \mapsto x_7$ | $x_8 \mapsto x_4$ | $x_8 \mapsto x_8$ |
| $x_9 \mapsto x_{10}$ | $x_9 \mapsto x_{11}$ | $x_9 \mapsto x_{12}$ | $x_9 \mapsto x_3$ | $x_9 \mapsto x_7$ | $x_9 \mapsto x_4$ | $x_9 \mapsto x_8$ | $x_9 \mapsto x_9$ |
| $x_{10} \mapsto x_{11}$ | $x_{10} \mapsto x_{12}$ | $x_{10} \mapsto x_3$ | $x_{10} \mapsto x_7$ | $x_{10} \mapsto x_4$ | $x_{10} \mapsto x_8$ | $x_{10} \mapsto x_9$ | $x_{10} \mapsto x_{10}$ |
| $x_{11} \mapsto x_{12}$ | $x_{11} \mapsto x_3$ | $x_{11} \mapsto x_7$ | $x_{11} \mapsto x_4$ | $x_{11} \mapsto x_8$ | $x_{11} \mapsto x_9$ | $x_{11} \mapsto x_{10}$ | $x_{11} \mapsto x_{11}$ |
| $x_{12} \mapsto x_3$ | $x_{12} \mapsto x_7$ | $x_{12} \mapsto x_4$ | $x_{12} \mapsto x_8$ | $x_{12} \mapsto x_9$ | $x_{12} \mapsto x_{10}$ | $x_{12} \mapsto x_{11}$ | $x_{12} \mapsto x_{12}$ |